\documentclass{birkjour}
\usepackage{amssymb,amsmath}
\usepackage{graphicx, wrapfig}

\begin{document}

\title[How Riemannian Manifolds Converge]{How Riemannian Manifolds Converge}

\author[Christina Sormani]{Christina Sormani\footnote{Research funded in part by PSC CUNY and NSF DMS 10060059}}
\address{%
Mathematics Department\\
CUNY Graduate Center\\
365 Fifth Avenue\\
NY NY 10014\\
USA\\ {\em and}\\
Mathematics Department\\
Lehman College\\
Bronx, NY 10468\\
USA}
\email{sormanic@member.ams.org}

\dedicatory{Dedicated to Jeff Cheeger for his 65th birthday}

\begin{abstract}
This is an intuitive survey of extrinsic and intrinsic notions
of convergence of manifolds complete with pictures of key examples
and a discussion of the properties associated with each notion.    
We begin with a description of three extrinsic
notions which have been applied to study sequences of 
submanifolds in Euclidean space: Hausdorff convergence of sets, 
flat convergence of integral currents, and weak convergence of
varifolds.  We next describe a variety of intrinsic notions
of convergence which have been applied to study sequences of 
compact Riemannian manifolds: Gromov-Hausdorff convergence
of metric spaces, convergence of metric measure spaces,
Instrinsic Flat convergence of integral current spaces, and ultralimits
of metric spaces.  We close with a speculative section addressing
possible notions of intrinsic varifold convergence, convergence of
Lorentzian manifolds and area convergence.
\end{abstract}

\maketitle


\newcommand {\Pind} {\hspace{.51cm}}  
\newcommand {\vare}{\varepsilon}
\newcommand {\grad}{\bigtriangledown}
\newcommand{\RR}{{\mathbb R}}
\newcommand{\SSS}{{\mathbb S}}
\newcommand{\EE}{{\mathbb E}}
\newcommand{\HH}{{\mathbb H}}
\newcommand{\ZZ}{{\mathbb Z}}
\newcommand{\NN}{{\mathbb N}}
\newcommand{\QQ}{{\mathbb Q}}
\newcommand{\be}{\begin{equation}}
\newcommand{\ee}{\end{equation}}
\newcommand{\Ricci}{\textrm{Ricci}}
\newcommand{\diam}{\textrm{diam}}
\newcommand{\mass}{\textrm{M}}
\newcommand{\vol}{\textrm{vol}}
\newcommand{\GHto}{\stackrel { \textrm{GH}}{\longrightarrow} }
\newcommand{\Fto}{\stackrel {\mathcal{F}}{\longrightarrow} }
\newcommand{\Hto}{\stackrel { \textrm{H}}{\longrightarrow} }
\newcommand{\spt}{\textrm{spt}}
\newcommand{\set}{\textrm{set}}


\newtheorem {theorem} {Theorem}
\newtheorem {proposition} [theorem] {Proposition}
\newtheorem {lemma} [theorem] {Lemma}
\newtheorem {definition} [theorem] {Definition}
\newtheorem {defn} [theorem] {Definition}
\newtheorem {corollary} [theorem] {Corollary}
\newtheorem {remark} [theorem] {Remark}
\newtheorem {note} [theorem] {Note}
\newtheorem {example} [theorem] {Example}
\newtheorem {conjecture} [theorem] {Conjecture}
\newtheorem {question} [theorem] {Question}


\section{Introduction}
The strong notions of smoothly or Lipschitz
converging manifolds have proven to be exceptionally useful when
studying manifolds with curvature and volume bounds, 
Einstein manifolds, 
isospectral manifolds of low dimensions, 
conformally equivalent manifolds, Ricci flow and
the Poincare conjecture, and even some questions in general relativity.  
{\em However
many open questions require weaker forms of convergence that do not
produce limit spaces that are manifolds themselves.}  Weaker notions
of convergence and new notions of limits have proven necessary in
the study of
manifolds with no curvature bounds or only lower bounds on Ricci or scalar 
curvature,
isospectral manifolds of higher dimension,
Ricci flow through singularities,
and general relativity.   
Here we survey a variety of weaker notions of  convergence and the corresponding
limit spaces covering
both well established concepts, newly discovered ones and speculations.

We begin with the
the convergence of submanifolds of Euclidean space as there is a wealth of 
different weak kinds of convergence that mathematicians have been applying 
for almost a century.  
Section~\ref{Sect-1.1} covers the {\em Hausdorff convergence of sets}, 
Section~\ref{Sect-1.2} covers Federer-Fleming's {\em flat convergence of integral currents}, and 
Section~\ref{Sect-1.3} covers Almgren's {\em weak convergence of varifolds}.  Each
has its own kind of limits and preserves different properties.  Each is useful for
exploring different kinds of problems.  While Hausdorff convergence is most well
known in the study of convex sets, flat convergence in the study minimal surfaces 
and varifold convergence in the study of mean curvature flow, each has appeared
in diverse applications.   Keep in mind that these are extrinsic notions of convergence
which depend upon how a manifold is located within the extrinsic space.

It is natural to believe that each should have
a corresponding intrinsic notion of convergence which should prove
useful for studying corresponding intrinsic questions about Riemannian
manifolds which do not lie in a common ambient space.  
The second section vaguely describes some of these
intrinsic notions of convergence with pictures illustrating key examples.
Section~\ref{Sect-2.1} covers  \em{Gromov-Hausdorff convergence of metric spaces}
which is an intrinsic Hausdorff distance,
Section~\ref{Sect-2.2} covers  various notions of \em{convergence
of metric measure spaces},
Section~\ref{Sect-2.3} covers the \em{intrinsic flat convergence on integral current spaces}, and
Section~\ref{Sect-2.4} covers the  weakest notion of all: \em{ultralimits of metric spaces}.  
Illustrated definitions and key examples will be given for each kind
of convergence.
This survey is not meant to provide a thorough rigorous definition of any of these
forms of convergence but rather to provide the flavor of each notion and direct
the reader to further resources.  

The survey closes with speculations on new notions of convergence:  
Section~\ref{Sect-3.1} describes difficulties arising when attempting to 
define an intrinsic varifold convergence.
Section~\ref{Sect-3.3}  describes the possible notion of area convergence and area spaces.  
Section~\ref{Sect-3.2} discusses the importance of weaker notions of convergence
for Lorentzian manifolds. 
The author attempts to include all key citations of initial work in these directions.


The author apologizes for the necessarily
incomplete bibliography and encourages the reader to consult mathscinet and
the arxiv for the most recent results in each area.  For some of the older
forms of convergence entire textbooks have been written focussing on
one application alone.  

\section{Converging Submanifolds of Euclidean Space}

There are numerous textbooks written covering the three notions of
extrinsic convergence we describe in this section.  Textbooks covering
all three notions include \cite{Morgan-text} and \cite{Lin-Yang-GMT-text}.
Morgan's text providing pictures and overviews of proofs with precise references
to theorems and proofs in Federer's classic text \cite{Federer}.  The focus
is on minimal surface theory.  
Another source with many pictures
and intuition regarding Plateau's problem, is Almgren's classic textbook
\cite{Almgren-Plateau}.
Lin-Yang's textbook \cite{Lin-Yang-GMT-text}
covers a wider variety of applications and provides a more modern
perspective incorporating recent work of Ambrosio-Kirchheim.  
An excellent resource on varifolds is Brakke's book \cite{Brakke}
which is freely available on his webpage.  Simon's classic
text \cite{Simon-text} is another indispensable resource.

\subsection{Hausdorff Convergence of Sets} \label{Sect-1.1}

The notion of Hausdorff distance dates back to the early 20th century.
Here we define it on an arbitrary metric space, $(Z, d_Z)$, although initially
it was defined on Euclidean space.
The Hausdorff distance between subsets  $A_1,A_2 \subset Z$ is 
 \be \label{eqn-Hausdorff-dist}
 d^Z_H(A_1, A_2) :=\inf\Big\{\, R \,{\bf :} 
 \, \, A_1 \subset T_R(A_2), \,\, A_2 \subset T_R(A_1) \Big\},
 \ee 
where $T_r(A):=\{y: \, \exists \, x\in A\,\,s.t.\,\, d_Z(x,y)<r\}$.
We write $A_j \Hto A$ iff $d^Z_H(A_j,A) \to 0$.

\begin{figure}[h] 
   \centering
   \includegraphics[width=4.6in]{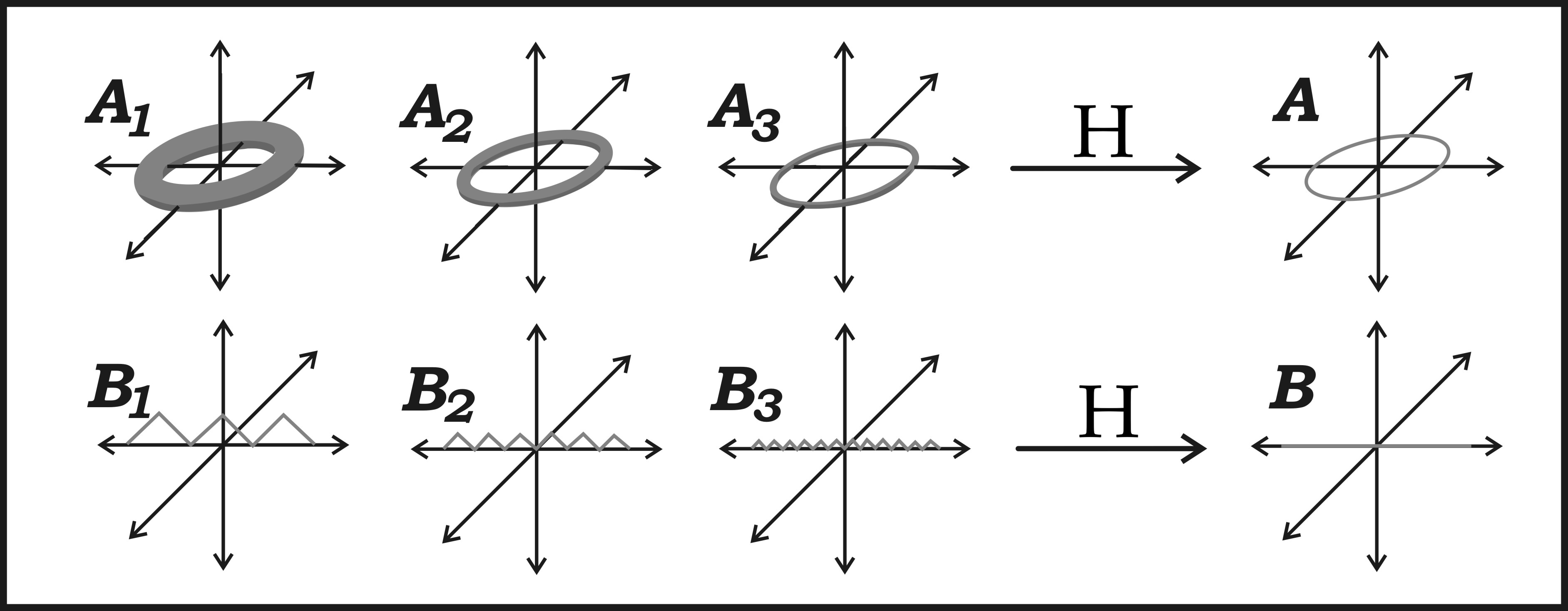} 
   \caption{Hausdorff Convergence}
   \label{fig-Hconv}
\end{figure}

In Figure~\ref{fig-Hconv}, we see
two famous sequences of submanifolds converging in the Hausdorff sense.  
The sequence $A_i$ are tori which converge to a circle $A$, depicting
how a sequence may lose both topology and dimension
in the limit.  Sequence $B_i$ of jagged paths in the $yz$ plane
converges to a straight segment, $B$, in the $y$ axis.  
Notice that $L(B_i)=\sqrt{2}$ while $L(B)=1$.  In addition to a sudden loss
of length under a limit, we lose all information about the derivative
(much like a $C_0$ limit).  In fact one may construct a sequence of jagged
curves that converge in the Hausdorff sense to the nowhere differentiable
Weierstrass function.

One property which is preserved by Hausdorff convergence is convexity.
For this reason, Hausdorff convergence has often been applied in the
study of convex subsets of Euclidean space.

In 1916, Blaschke proved that if a sequence of nonempty compact sets
$K_i$ lie in a ball in Euclidean space, then a subsequence converges in
the Hausdorff sense to a nonempty compact set \cite{Blaschke-1916}.  If
the sets are connected, then their limit is connected.  Path connectedness,
however, is not preserved as can be seen by a sequence approaching
the well known set: 
\be
\{(x,sin(1/x)): \, x\in (0,1]\}\cap \{(0,y): \, y\in [-1,1]\}.
\ee

Dimension and measure are also not well controlled.
The $k$ dimensional Hausdorff measure of a set $X$
is defined by covering $X$ with countable collections of
sets $C_i$ of small diameter: 
\begin{equation*} 
\mathcal{H}_k(X):=\lim_{r\to 0} \inf\left\{\sum_{i=1}^\infty \alpha(k) 
\left(\frac{\diam(C_i)}{2}\right)^k \,:\,\, X\subset \bigcup_{i=1}^\infty C_i, \,\, \diam(C_i)<r\right\}
\end{equation*}
where $\alpha(k)$ is the volume of a unit ball of dimension $k$ in
Euclidean space.  For a beautiful exposition of this notion see \cite{Morgan-text}.
When $X$ is a submanifold of dimension $k$, then 
$\mathcal{H}_k(X)$ is just the $k$ dimensional Lebesgue measure: when $k=1$
it is the length, when $k=2$ it is the area and so on.  

\begin{figure}[h] 
   \centering
   \includegraphics[width=4.6in]{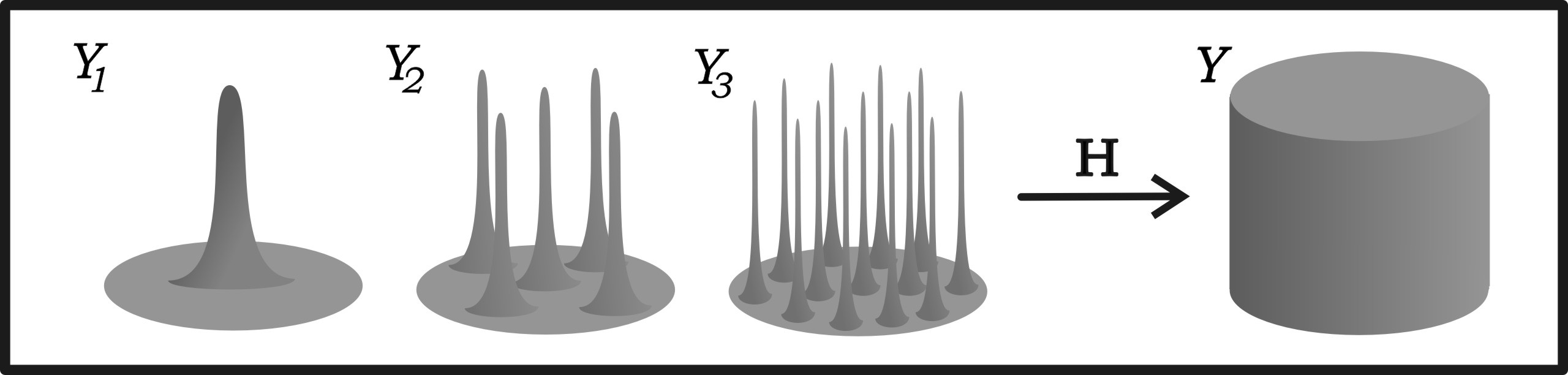} 
\caption{Disk of Many Splines}
   \label{fig-splines}
\end{figure}

The Hausdorff dimension of a set, $X$, is 
\be \label{eqn-def-Hausdorff-dimension}
H_{\textrm{dim}}(X):=\inf\{k \in(0,\infty): \mathcal{H}_k(X)=0\}.
\ee
A compact k-dimensional submanifold, $M^k$, has $H_{\textrm{dim}}(M)=k$.
Notice that in Figure~\ref{fig-Hconv}, $H_{\textrm{dim}}(A_j)=2$ while
$H_{\textrm{dim}}(A)=1$.  It is also possible for the dimension to
go up in the limit.  In Figure~\ref{fig-splines}, $H_{\textrm{dim}}(Y_j)=2$ while
$H_{\textrm{dim}}(Y)=3$.  It is also possible for a sequence of $1$
dimensional submanifolds, to converge to a space of fractional
dimension like the von Koch curve. 
Such a sequence of curves must have length diverging to infinity.

In fact, if $X_j$ have uniformly bounded
length, $\mathcal{H}_1(X_j) < C$, and $X_j$ converge to $X$ in the
Hausdorff sense, then 
\be \label{eqn-H-length}
\mathcal{H}_1(X) \le \liminf_{j\to\infty}\mathcal{H}_1(X_j).
\ee
Note that the sequence $B_j$ depicted in Figure~\ref{fig-Hconv} has
$\mathcal{H}_1(B_j)=\sqrt{2}$ for all $j$, but $\mathcal{H}_1(B)=1$,
so lower semicontinuity is the best one can do.
This is not true for higher dimensions as can be seen in Figure~\ref{fig-splines},
where $\lim_{j\to\infty} \mathcal{H}_2(Y_j)=1$ yet $\mathcal{H}_2(Y)=\infty$.

Figure~\ref{fig-splines} depicts a famous example of ``a disk with splines" described in Almgren's text on Plateau's problem \cite{Almgren-Plateau}.  In Plateau's problem one is given a closed curve and is asked to find the surface of smallest area spanning that curve.
In this example, the sequence of smooth surfaces, $Y_j$,
have a common boundary, $\partial Y_j=\{(x,y,0):\, x^2+y^2=1\}$ 
which is a given closed curve, and the area of the $Y_j$ is approaching the minimal area
filling that circle.  However the sequence does not converge in the Hausdorff
sense to the flat disk, $D^2$, which is the solution to the Plateau problem for a circle.  
Instead, due to the many splines, the sequence converges in the Hausdorff
sense to the solid cylinder.  

The disk with many splines example convinced mathematicians that the Hausdorff
distance was not useful in the study of minimal surfaces.  New notions of
convergence had to be defined. 


\subsection{Flat Convergence of Integral Currents} \label{Sect-1.2}

Federer and Fleming introduced the notion of an integral current and
the flat convergence of integral currents to deal with examples
like the disk with many splines depicted in Figure~\ref{fig-splines}.
When viewed as integral currents, the submanifolds depicted in
that figure converge in the flat sense to the disk.  All the splines
disappear in the limit.

A current, $T$, is a linear functional on smooth $k$ forms.
Any compact oriented submanifold, $M$, of dimension $k$ with a smooth
compact boundary, may be viewed an a $k$ dimensional 
current, $T$, defined by 
\be
T(\omega):=\int_M  \omega.
\ee
Notice that in Figure~\ref{fig-splines}, we have 
\be
\lim_{j\to\infty}\int_{Y_j}\omega = \int_{D^2} \omega
\ee
for any smooth differentiable $3$ form, $\omega$.  So, viewed as currents,
$Y_j$ converge weakly to the flat disk $D^2$.

Federer-Fleming proved that any sequence of compact $k$ dimensional oriented
submanifolds, $M_j$, in a disk in Euclidean space, with a uniform
upper bound on $\mathcal{H}^k(M_j)\le V_0$ and a uniform
upper bound on $\mathcal{H}^{k-1}(\partial M_j) \le A_0$, has a subsequence
which converges when viewed as currents in the weak sense.
The limit is an ``integral current".

An integral current, $T$, is a current with a canonical set, $R$, and a
multiplicity function $\theta$.  The canonical set $R$ is countably 
$\mathcal{H}^k$ rectifiable,
which means it is contained in the image of a countable collection of Lipschitz
maps, $\varphi_i:E_i\to R$ from Borel subsets, $E_i$, of $k$ dimensional Euclidean space.  The multiplicity
function $\theta$ is an integer valued Borel function.  We define
\be \label{eqn-def-integral-current}
T(\omega) \,:= \,\int_R \theta \, \omega \,=\, \sum_{i=1}^\infty \,\int_{E_i} \theta\circ\varphi_i \,\,\,
\varphi_i^*\omega.
\ee
It is further required that an integral current have finite mass
\be \label{eqn-def-mass}
\mass(T) := \int_R \theta \, d\mathcal{H}^k
\ee
and that the boundary, $\partial T$, defined by
$
\partial T (\omega) := T( d\omega)
$
have finite mass, $\mass (\partial T)<\infty$.  Note that Federer-Fleming 
and Ambrosio-Kirchheim proved
this implies that
$\partial T$ also has a rectifiable canonical set although it is
one dimension lower than the dimension of $T$.   \footnote{See \cite{Morgan-text} and \cite{Lin-Yang-GMT-text} for more details.
Note there are slight differences in the definition as \cite{Morgan-text} follows
Federer-Fleming \cite{FF}, while \cite{Lin-Yang-GMT-text} follows the newer
version introduced by Ambrosio-Kirchheim \cite{AK}.}

When a submanifold
$M$ is viewed as an integral current $T$, then $M$ itself is the
canonical set, it has multiplicity $1$, the boundary, $\partial T$,
is just $\partial M$ viewed as an integral current and the mass,
$\mass(T)$, is just the volume of $M$.

Also included as a $k$ dimensional integral current is the
$0$ current.  In Figure~\ref{fig-Hconv}, the sequence $A_j$ may be viewed as integral
currents.  They converge in the weak sense as integral currents to
the ${ 0}$ current.  More generally, if $\mathcal{H}^k(M_j)$ decreases to
zero, then the weak limit of the $k$ dimensional submanifolds $M_j$
viewed as integral currents is also the $0$ current.   In fact, whenever
$M_j$ have a uniform upper bound on $\mathcal{H}^k(M_j)$ and 
$M_j=\partial N_j$ where $\mathcal{H}^{k+1}(N_j) \to 0$, then
$M_j$ viewed as integral currents also converge weakly to the $0$ current.
The sequence $A_j$ depicted in Figure~\ref{fig-Hconv} are the
boundaries of solid tori whose volumes decrease to $0$.
It was not actually necessary that their areas decrease to $0$.
This idea of filling in the manifold to assess where it converges makes
it much easier to see the limits of integral currents and leads naturally
to the following definition.

The flat distance between two $k$ dimensional integral currents, 
$T_1$ and $T_2$, is defined by
  \be \label{eqn-def-flat}
 d_{\mathcal{F}}(T_1,T_2) =\inf \Big\{\mass(A) + \mass(B)
 \,{\bf {:}} \,\, \textrm{\em{int curr }}A,B \,\,s.t.\,\,  A +\partial B=T_1-T_2\Big\},
 \ee
where the infimum is taken over all $k$ dimensional integral currents, $A$,
and all $k+1$ dimensional integral currents, $B$.  \footnote{See prior footnote.}

\begin{wrapfigure}{r}{2.2in} 
   \begin{center}
   \includegraphics[width=2in]{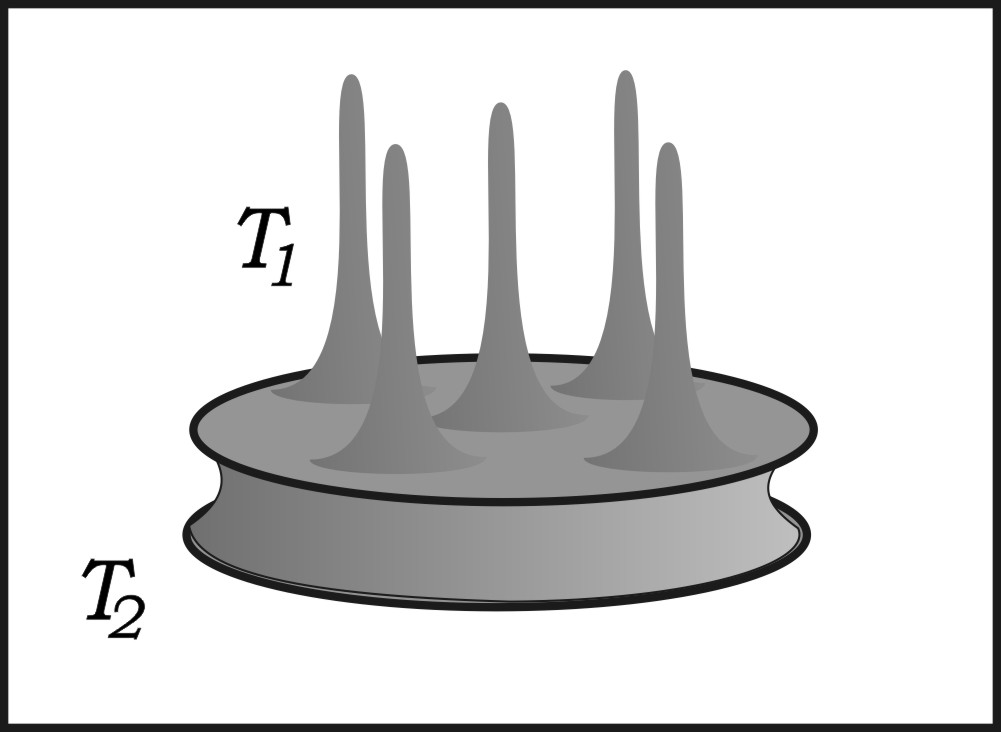} 
\end{center}
\end{wrapfigure}

In the figure to the right 
we see a choice of $A$ with small area
that looks like a catenoid, and then we fill in the space between the flat disk
$T_2$ and the disk with many splines, $T_1$, to define $B$.  If the 
surface with the splines actually shares the same boundary with the disk
(as it does in Figure~\ref{fig-splines}), then we can take $A$ to be
the $0$ current and $B$ just to be a filling (the region in between them).
It is easy to see that the sequence of manifolds $Y_j$ depicted in
Figure~\ref{fig-splines} converges to the disk $D^2$ using the flat
distance. 

In Figure~\ref{fig-Hconv} the sequence $B_j$ viewed as $1$ dimensional
integral currents, $T_j$, converges to the Hausdorff limit $B$ viewed as a one
dimensional current, $T$ (as long as they are given the same orientation
left to right).  Here $T_j-T$ is again a cycle and we can find surfaces 
viewed as $2$ dimensional currents, $S_j$,
such that $\partial S_j=T_j-T$ whose areas, $\mass(S_j)$, converge to $0$.

Federer-Fleming proved that when a sequence of integral currents
has a uniform upper bound on mass and on the mass of their boundaries,
then they converge weakly iff they converge with respect to the flat distance.
Their compactness theorem can now be restated as follows: {\em if $T_j$
is a sequence of $k$ dimensional integral currents supported in a compact subset
with $\mass(T_j)\le V_0$
and $\mass(\partial T_j)\le A_0$ then a subsequence converges in the
flat sense to an integral current space.}

One of the beauties of this theorem is that the limit space is
rectifiable with finite mass,
and in fact the mass is lower semicontinuous.  This makes flat convergence
an ideal notion when studying Plateau's problem.  The limit integral current
has the same boundary as the sequence and minimal area.  

\begin{figure}[h] 
   \centering
   \includegraphics[width=4.6in]{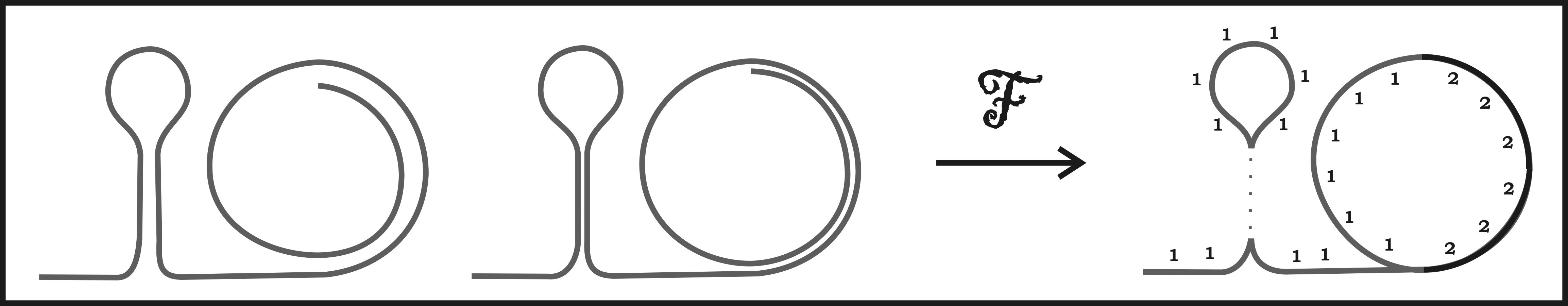} 
   \caption{Doubling and Cancellation of Flat Limits}
   \label{fig-cancels}
\end{figure}

In addition integral currents have a notion of an approximate tangent plane,
which exists almost everywhere and is a subspace of the same
dimension as the current.  When $M$ is a submanifold, this approximate
tangent plane is the usual tangent plane and exists everywhere inside
$M$.  In 1966 Almgren was able to apply the strong control on the approximate
tangent planes to obtain even stronger regularity results for the limits
achieved when working on Plateau's problem \cite{Almgren-1966}.  
The limit space ends up having multiplicity $1$.  In general the limit
space may have higher multiplicity or regions with higher multiplicity
[Figure~\ref{fig-cancels}], which we will discuss in more detail later.

In 1999, Ambrosio-Kirchheim extended the notion of an integral current
from Euclidean space to arbitrary metric spaces, proving the
existence of a solution to Plateau's problem on Banach spaces \cite{AK}.
A key difficulty was that on a metric space, $Z$, there  is no notion
of a differential form.   They applied a notion of DeGiorgi \cite{DeGiorgi}, replacing
a $k$ dimensional differential form like
$\omega=fdx_1\wedge dx_2\wedge\cdots\wedge dx_k$ with a $k+1$
tuple, $\omega=(f,x_1, x_2, ...x_k)$, of Lipschitz functions
satisfying a few rules, including
$
d\omega:\,\,=\,\,(1,f,x_1,x_2,...x_k).$

An integral current, $T$, is defined as a linear functional on $k+1$ tuples
using a rectifiable set, $R$, and a multiplicity function, $\theta$, so that
as in (\ref{eqn-def-integral-current}), we have
\be
T(\omega) \,\,: = \,\,\sum_{i=1}^\infty \int_{E_i} \theta\circ\varphi_i \,\,\,
\varphi_i^*\omega.
\ee
where $\varphi_i^*\omega \,:= \,f\circ \varphi \,\,d(x_1\circ \varphi)\wedge d(x_2\circ\varphi)
\wedge \cdots \wedge d(x_k\circ \varphi)$ is defined almost everywhere
by Rademacher's Theorem.
They can then define mass exactly as in (\ref{eqn-def-mass}) and
boundary exactly as before as well.   They require integral currents
to have bounded mass and their boundaries to have bounded mass.\footnote{See \cite{AK} or \cite{SorWen2} for the precise definition.}

Ambrosio-Kirchheim prove that if $T_j$ are integral currents in a compact metric 
space and $\mass(T_j) \le V_0$ and $\mass(\partial T_j)\le A_0$ then a
subsequence converges in the weak sense to an integral current.  Furthermore
the mass is lower semicontinuous.  There is also a notion of an approximate
tangent plane, however, here the approximate tangent plane is a normed
space.  The norm is defined using the metric differential described in earlier work by
Korevaar-Schoen on harmonic maps and also by Kirchheim
on rectifiable space \cite{Korevaar-Schoen}\cite{Kirchheim}. 

Wenger extended the class of metric spaces $Z$ which have a solution
to Plateau's problem and  defined a flat distance exactly as in (\ref{eqn-def-flat}).
We will use the notation $d_{\mathcal{F}}^Z$ for the flat distance in $Z$ to be
consistent with $d_H^Z$ denoting the Hausdorff distance in $Z$.
He proved that on this larger class of spaces, which includes Banach spaces,
when $\mass(T_j) \le V_0$ and $\mass(\partial T_j)\le A_0$ then 
$T_j$ converge weakly to $T$ iff $T_j$ converges in the flat sense to $T$
\cite{Wenger-Plateau}\cite{Wenger-flat}.
Recall that by Rademacher's Theorem, any separable metric space, $Z$,
isometrically embeds into a Banach space, $\textrm{W}$, so we can always study
the convergence of integral currents in $Z$ using the flat distance
of the push forwards of the currents in Banach space, $d_{\mathcal{F}}^{\textrm{W}}$.
If $\varphi:Z \to \textrm{W}$, then the push forward of a current $T$ in $Z$, is a current
$\varphi_\#T$ in $\textrm{W}$ defined by
$
\varphi_\# T \, (f,x_1,...x_k) \,\, := \,\, T \, (f\circ \varphi, x_1\circ\varphi, ... x_k \circ \varphi).
$

Another extension of the notion of integral currents is that of the
$\ZZ_n$ integral currents where the mutiplicity function takes values
in $\ZZ_n$.  This was introduced by Fleming \cite{Fleming-finite}
and has been extended to arbitrary metric spaces by Ambrosio and
Wenger \cite{Ambrosio-Wenger}.  One application for this notion
is the study of minimal graphs (c.f. \cite{Almgren-Plateau}).
To understand this better we will go over an example using Federer-Fleming's notation.

In Figure~\ref{fig-cancels}, the sequence of embedded curves $C_i$
converges in $C^1$ to a limit curve $C$ which is not embedded. 
Viewed as $1$ dimensional currents, $T_i$, they converge in the flat sense to $T$.  
It is perhaps easier to understand Figure~\ref{fig-cancels} using weak
convergence rather than flat convergence.  Note that
\begin{equation*}
T_i(\,f_1(x,y)dx+ f_2(x,y) dy\,)= \int_0^1 (f_1\circ C_i)(t)\, 
d(x\circ C_i) + (f_2\circ C_i)(t)\,  d(y\circ C_i)
\end{equation*}
converges to 
\begin{equation*}
T(\,f_1(x,y)dx+ f_2(x,y) dy\,)\,\,= \,\,\int_0^1 (f_1\circ C)(t)\, 
d(x\circ C) + (f_2\circ C)(t)\,  d(y\circ C)
\end{equation*}
Notice that $T$ has multiplicity $1$
where $C$ does not overlap itself. 
It has multiplicity $2$ on the segment $C$ covers twice in the same
direction.  The segment depicted with dashes where $C$ passes
in opposite directions is not part of the canonical set of $T$ because
the integration cancels there.    
 If one uses $\ZZ_n$ valued coefficients
we get the same limit space for $n>2$, but for $n=2$ the doubled up
segment disappears as well as the cancelled dashed segment.

The cancellation depicted in Figure~\ref{fig-cancels} is problematic for
some applications.  Note that the same effect can occur in higher
dimensions, by extending the curve to a sheet in $\RR^3$.
While often it is useful for thin splines to disappear
in the hopes of preserving dimension, canceling sheets which overlap
with opposing orientations is not necessary to obtain a rectifiable limit
space.  This is seen using the notion of a varifold.

\subsection{Weak Convergence of Varifolds}  \label{Sect-1.3}

Varifolds were introduced by Almgren in 1964 \cite{Almgren-notes}.  
Significant further work was completed by Allard in 1972 \cite{Allard-1972}.
Almgren's goal was to define
a new notion of convergence for submanifolds which had rectifiable limit
spaces and a notion of tangent planes for those limit spaces, but did
not have the kind of cancellation that occurs when taking flat limits of
integral currents.  Under varifold convergence, the sequence of
$1$ dimensional manifolds depicted in Figure~\ref{fig-cancels}, converges
to the rectifiable set with weight $1$ everywhere including the dashed 
segment.  There is no notion of orientation on the limit, but there
are tangent planes almost everywhere.

A $k$ dimensional varifold is a Radon measure on $\RR^N \times \Gamma(k,N)$ where
$\Gamma(k,N)$ is the space of $k$ subspaces of Euclidean space, $\RR^N$.
\footnote{In some books $\Gamma(k,N)$ is written as $GL(k,N)$}
A submanifold $M^k\subset \RR^n$, may be viewed as the varifold, $V$,
defined on any $W\subset \RR^N \times \Gamma(k,N)$ as follows:
\be \label{eqn-def-varifold}
V(W) := \mathcal{H}^k \left( \, W\,  \cap\, \left\{ (x, T_xM): \, x\in M \subset \RR^m\right\} \right),
\ee
where $T_xM$ is the tangent space to $x$ translated to the origin. 

A sequence of varifolds, $V_j$, is said to converge weakly to a varifold, $V$,
\be \label{eqn-def-varifold-convergence}
\int f \, dV_j\,\, \to \,\,\int f \, dV \qquad \forall \, f \, \in \, C_0(\RR^N\times \Gamma(k,N)).
\ee

For example, let us examine the sequence of curves $C_j$ which converge $C^1$
to the curve $C$ in Figure~\ref{fig-cancels}.
They may be viewed as integral currents $V_j$,
\be \label{eqn-ex-varifold-1}
V_j(W)=\mathcal{H}^1\left( W \, \cap \,
\left\{ \left( C_j(t), \pm C_j'(t)/|C_j'(t)| \right): \, t\in [0,1]\right\} \, \right),
\ee
where we view points in $\Gamma(1,2)$ as $\pm v$ where $v$ in a unit
$2$ vector.  Then
\be \label{eqn-ex-varifold-2}
\int f(x,v) dV_j \, = \, \int_0^1 f(C_j(t), \pm C_j'(t)/|C_j'(t)|) \, dt
\ee
converges to $\int_0^1 f(C(t), \pm C(t)/|C(t)|) dt =\int f(x.v)\, dV$
for the varifold $V$ defined by
\be \label{eqn-ex-varifold-3}
V(W) := \int _0^1 \chi_W(C(t), \pm C'(t)/|C'(t)|) \,dt
\ee
where $\chi_W$ is the indicator function of $W$.  This varifold corresponds
to viewing the limit curve $C$ as having weight $1$ on the segments where
it doesn't overlap itself and weight $2$ on the segments where it does overlap.
There is no cancellation.  The dashed segment has weight $2$.

Not all limit varifolds end up with tangent planes that align well with the
rectifiable set as they do in (\ref{eqn-ex-varifold-3}).  If we examine instead
the sequence $B_j$ in Figure~\ref{fig-Hconv}, as a sequence of curves $C_j$
converging to a curve $C$, we see that the limit varifolds, $V$, has the form
\be \label{eqn-ex-varifold-4}
V(W) := \,\,a \int _0^1 \chi_W(C(t), \pm av_1) \,dt \,\,+\,\,a \int _0^1 \chi_W(C(t), \pm av_2) \,dt
\ee
where $v_1=(1,1)$, $v_2=(-1,1)$ and $a =\sqrt{2}/2$ because $C_i'(t)$ is always
in the direction of $v_1$ or $v_2$. 
This effect was observed by Young
and is discussed at length in \cite{Lin-Yang-GMT-text}.   Varifolds like
$V$ are of importance but are not considered to be integral varifolds
because $C'(t)$ is unrelated to $v_1$ and $v_2$.
 
An integral varifold, $V$, is defined as a positive integer weighted countable sum
of varifolds, $V_j$, which are
defined by embedded submanifolds $M_j$ as in (\ref{eqn-def-varifold}).
While $V$ is a measure on $\RR^N\times \Gamma(N,k)$, it
defines a natural measure $||V||$ on $\RR^N$,
\be \label{eqn-def-varifold-density}
||V||(A) := V(A \times \Gamma(k,N)).
\ee 
When $V$ is defined by a submanifold, $M$, then $||V||(A)=\mathcal{H}^k(A \cap M)$.

A $k$ dimensional varifold $V$  is said to have a tangent plane, $T\in \Gamma(k,N)$ with
multiplicity $\theta\in (0,\infty)$ at a point  $x\in \RR^N$ if a sequence of rescalings
of $V$ about the point $x$  converges to $\theta(x) T$. \footnote{See \cite{Lin-Yang-GMT-text}
for a more precise definition.}  
An integral varifold, $V$, 
has a tangent space $||V||$ almost everywhere on $\RR^N$.

Varifolds do not have a notion of boundary.  Instead a varifold, $V$, has a notion of
first variation, $\delta V$, which is a functional that maps functions 
$f\in C_0^\infty(\RR^N,\RR^N)$
to $\RR$.  There is no room here to set up the background for
a precise definition \footnote{See \cite{Lin-Yang-GMT-text} 6.2 for the full definition.}, 
so we just describe
$\delta V$ when $V$  is defined by a $C^2$ submanifold $M$
with boundary $\partial M$.   In this case, 
 $\delta V (f)$  is the first variation in the area of $M$ as it is flowed through
 a diffeomorphism defined by the vector field $f$:
 \be \label{eqn-def-first-var}
 \delta V(f) \, = \, - \int_M \, f(x) \cdot H(x) \, d\mathcal{H}^k(x) \,\,+\,\, \int_{\partial M} \, f(x)\cdot \eta(x)\, 
 d\mathcal{H}^{k-1}(x)
 \ee
 where $H(x)$ is the mean curvature at $x\in M$ and $\eta(x)$ is the outward normal
 at $x \in \partial M$.
 
 A varifold, $V$, is said to be stationary when $\delta V=0$.       
 
 A varifold has bounded first variation if
 \be
 ||\delta V|| (W) := \sup \{\delta V(f): \, f\in C_0^\infty(\RR^N, \RR^n), \, |f|\le 1,
 \, \textrm{spt}(f) \subset A \} <\infty.
 \ee
 When $V$ corresponds to a submanifold, $M$, as above, then
 \be
 ||\delta V|| (W) \, = \, \mathcal{H}^{k-1}(W \cap \partial M) \, + \,
 \int_{W \cap M} |H(x)| \, d\mathcal{H}^m(x). 
 \ee
 There is an isoperimetric inequality for varifolds.
  
 Allard proved that if a sequence of integral varifolds, $V_j$, with a uniform 
 bound on  $||\delta V_j||(W)$ depending only on $W$, converges
 weakly to $V$, then $V$ is an integral varifold as well \cite{Allard-1972}.
 In particular a weakly converging sequence of minimal surfaces, $M_j$, with a uniform
 upper bound on the length of their boundaries, $\partial M_j$, 
 converges to an integral varifold.  For more general submanifolds, one
 needs only uniformly control the volumes of the boundaries and the $L^1$ norms
 of the mean curvatures.
   
One key advantage of varifolds is that they have a notion of mean curvature 
defined using $\delta V$
 and an integral similar to the one in (\ref{eqn-def-first-var}).  This notion
 is used to define the Brakke flow, a mean curvature flow past singularities \cite{Brakke}.

Recently White has set up a natural map, $F$, from integral varifolds to $\ZZ_2$
flat chains, which basically preserves the rectifiable set and takes the integer
valued weight to a $\ZZ_2$ weight.  Each $\ZZ_2$ flat chain corresponds
uniquely to a $\ZZ_2$ integral current.
He has proven that if a sequence of submanifolds,
$M_j$
viewed as varifolds converge weakly to an integral varifold $V$
and satisfy the conditions of Allard's compactness theorem,
and if further $\partial M_j$ converge as $\ZZ_2$ integral currents, 
then $M_j$ viewed as $\ZZ_2$ integral currents converge to a 
$\ZZ_2$ integral current $T$ corresponding to $F(V)$.
 This is easily seen to be the situation for
the sequence depicted in Figure~\ref{fig-cancels}.  \cite{White-2009}

\section{Intrinsic Convergence of Riemannian Manifolds}

When studying sequences of Riemannian manifolds, $M_j$, the strongest notions
of convergence require that $M_j$ be diffeomorphic to the limit space $M$
with the metrics, $g_j$, converging smoothly:
\be
\exists \varphi_j: M \to M_j \textrm{ such that } \varphi_{j}^* g_j \to g.
\ee
  In this survey we are concerned with 
sequences of Riemannian manifolds which do not converge in such a strong sense.  
Here we describe a few weaker notions of convergence which allow us
to better understand sequences which do not converge strongly.  We begin 
with a pair of motivating examples.

The sequence of flat tori, $M_j=S^1_\pi \times S^1_{\pi/j}$, has volume converging
to $0$: $\vol(M_j)=2\pi(2\pi/j) \to 0$.  
Sequences with this property are called collapsing sequences.  They
do not converge in a strong sense to a limit which is also a torus.  Intuitively one
would hope to define a weaker notion of convergence in which these tori
converge to a circle.   Examples like this lead to Gromov's
notion of an intrinsic Hausdorff  convergence. 

If one views Figure~\ref{fig-splines} as a sequence of Riemannian disks
with splines, $M_j$, with the induced Riemannian metrics, they do not converge 
in a strong sense to a flat Riemannian disk.  While they are diffeomorphic to
the disk and the volumes are converging, $\vol(M_j)\to \vol(D^2)$,
the metrics do not converge smoothly.  Examples like this lead to the notion
of intrinsic flat convergence.

In this section we first present Gromov-Hausdorff convergence, then
metric measure convergence, then the intrinsic flat convergence and finally,
weakest of all, the notion of an ultralimit.  We include a few key examples,
applications and further resources for each notion.

%
%
%
%
\subsection{Gromov-Hausdorff Convergence of Metric Spaces}\label{Sect-2.1}

In 1981,  Gromov
introduced an intrinsic Hausdorff convergence for sequences of
metric spaces \cite{Gromov-1981}.  A few excellent references are
Gromov's book \cite{Gromov-metric},
 the textbook of Burago-Burago-Ivanov\cite{BBI},
 Fukaya's survey \cite{Fukaya-survey} and Bridson-Haefliger's
 book \cite{BH}.

\begin{figure}[h] 
   \centering
   \includegraphics[width=4.6in]{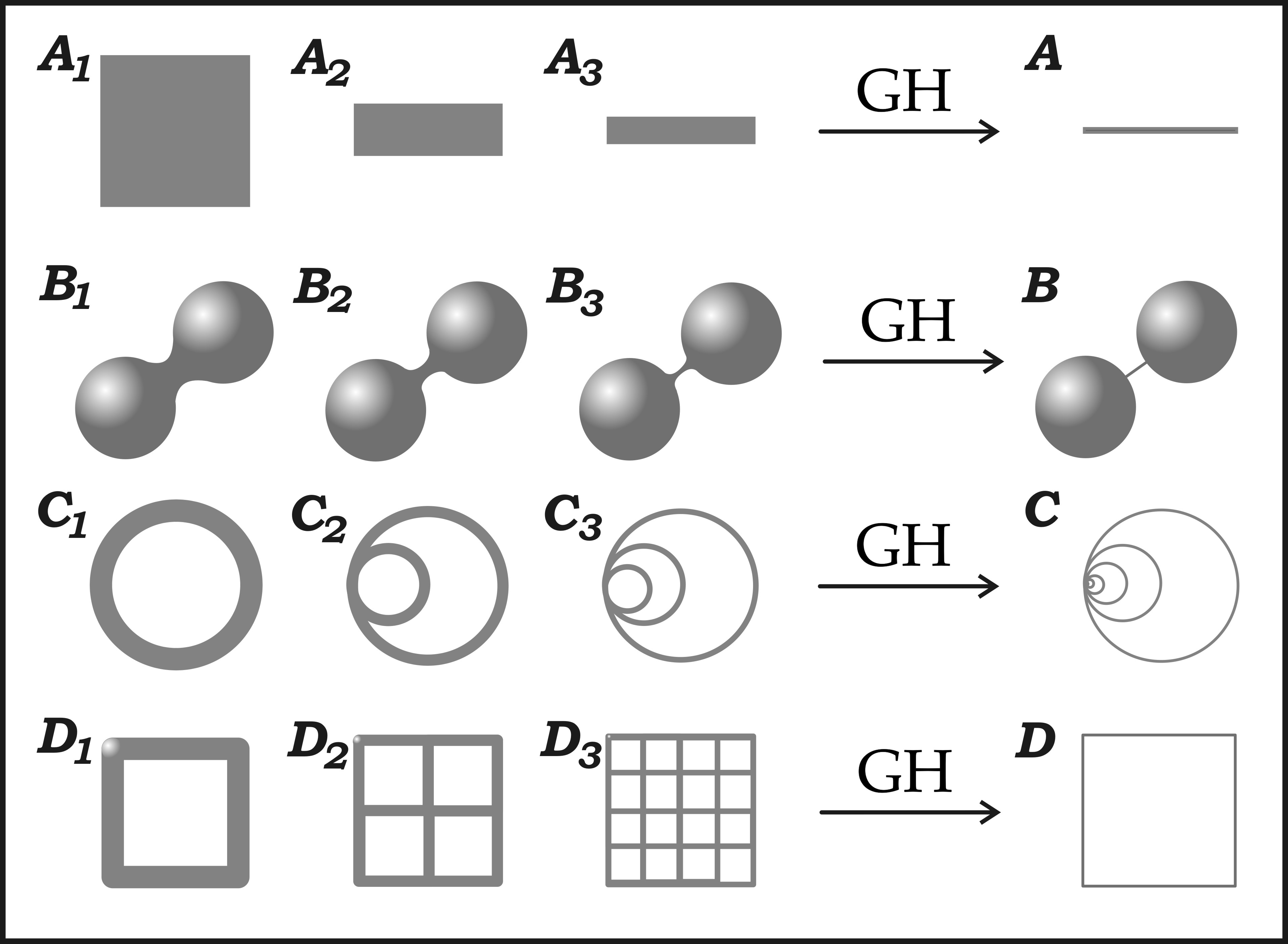} 
   \caption{Gromov-Hausdorff Convergence}
      \label{fig-GHconv}
\end{figure}

The Gromov-Hausdorff distance is 
 defined between any pair of compact metric spaces, 
  \be \label{eqn-GH-dist}
 d_{GH}(M_1, M_2)= \inf \Big\{ d^Z_H\left(\varphi_1(M_1), \varphi_2(M_2)\right)
 \,{\bf:} \,\,
 isom\,\, \varphi_i: M_i \to Z\Big\}
 \ee
where the infimum is taken over all metric spaces, $Z$, and all
isometric embeddings, $\varphi_i:M_i \to Z$.  
An isometric embedding, $\varphi:X \to Z$
satisfies
\be \label{eqn-isom-embed}
d_Z(\varphi(x_1), \varphi(x_2)) = d_X(x_1, x_2) \qquad \forall x_1, x_2 \in X.
\ee
We write $M_i \GHto X$ iff $d_{GH}(M_j, X) \to 0$.
See Figure~\ref{fig-GHconv}.

 The sequences of Riemannian manifolds depicted in Figure~\ref{fig-GHconv}
 reveal a variety of properties that are not conserved under Gromov-Hausdorff
 convergence.   The first sequence $A_j$ are the flat tori $S^1_\pi\times S^1_{\pi/j}$ 
 converging to a circle $A=S^1_\pi$.  To see this one takes the common space 
 $Z_j=A_j$ and isometrically embeds $A$ into $Z_j$ so that 
 \be
 d_{GH}(A_j, A) \le d_H^{Z_j}(A_j, \varphi_j(A))= \pi/(2j) \to 0.
 \ee
 Here we see the topology and Hausdorff dimension may decrease in the limit.
 
Note that if $M_1$ and $M_2$ are compact then $d_{GH}(M_1, M_2)=0$
iff $M_1$ and $M_2$ are isometric.  The Gromov-Hausdorff distance
between $M_1$ and $M_2$ is almost $0$ iff there is an almost isometry 
$f:M_1 \to M_2$ satisfying
\be \label{eqn-def-almost-isometry}
|d_1(x,y) - d_2(f(x),f(y))|<\epsilon \textrm{ and }
M_2\subset T_\epsilon(f(M_1)).
\ee
Note that an almost isometry need not be continuous.
In Figure~\ref{fig-GHconv}, it is easy to construct $\epsilon_j$
almost isometries, $\varphi_{j}: B \to B_j$, such that $\epsilon_j \to 0$
and conclude that $B_j \GHto B$.
This example 
 reveals that the Gromov-Hausdorff
 limit of a sequence of Riemannian manifolds may not be a Riemannian manifold.
 
 Gromov-Hausdorff limits of Riemannian manifolds are geodesic metric spaces.
 This means that the distance between any pair of points is equal to the length
 of the shortest curve between them.  The shortest curve exists and is called
 a minimal geodesic.  As in Riemannian geometry, a curve $\gamma$ is called
 a geodesic if for every $t$, there is an $\epsilon>0$ such that $\gamma$
 restricted to $[t-\epsilon, t+\epsilon]$ is a minimal geodesic.\footnote{Alexandrov 
 space geometers use the term ``geodesic" to refer to 
 a ``minimal geodesic" and ``local geodesic" to refer to a geodesic.}
 
 In the third sequence of Figure~\ref{fig-GHconv},  
 $C_j$ are also the boundaries of increasingly thin
 tubular neighborhoods.  The limit, $D$, is the Hawaii Ring, a metric space of
 infinite topological type that has no universal cover.  For more about the
 topology of Gromov-Hausdorff converging sequences of manifolds
 see \cite{SorWei1} \cite{SorWei3} and \cite{ShnSor2}. 
 
  In the last sequence of Figure~\ref{fig-GHconv}, 
 $D_j$ are the smoothed boundaries of increasing
 thin tubular neighborhoods of increasingly dense grids.  They converge to a 
 square $D=[0,1]\times[0,1]$.  
 The metric on the square is the taxicab metric (also called the
 $l_1$ metric):
 \be
 d_C((x_1,y_1),(x_2,y_2)) := |x_1-x_2| + |y_1-y_2|.
 \ee
 This is easiest to see by showing the $D_j$ are Gromov-Hausdorff
 close to their grids and that the grids converge to the taxicab square.
   
 Gromov proved that sequences of Riemannian manifolds, $M_j$, with uniform
 upper bounds on their diameter and on the number, $N(r)$, of disjoint balls 
 of radius $r$ have subsequences which converge in the Gromov-Hausdorff
 sense to a compact geodesic space, $Y$ \cite{Gromov-1981}.  
 Conversely, when $M_j$ converge to a compact $Y$, $N(r)$ is uniformly
 bounded.  Consequently, if one views the sequence of disks with splines 
 depicted in Figure~\ref{fig-splines} as Riemannian manifolds (with the
 intrinsic distance), the sequence
 does not converge in the Gromov-Hausdorff sense: a ball of radius $1/2$
 about the tip of a spline does not intersect with the a ball of radius $1/2$
 about the tip of another spline.  As the number of splines approaches infinity
 so does the number of disjoint balls of radius $1/2$. 
 
 By the Bishop-Gromov volume
 comparison theorem, sequences of manifolds, $M_j$, with uniform lower bounds
 on Ricci curvature and upper bounds on diameter satisfy these compactness criteria 
 \cite{Gromov-1981}.  This includes, for example, the sequence of flat tori, $A_j$,
 in Figure~\ref{fig-GHconv} as well as any other collapsing sequence of
 manifolds with bounded sectional curvature.  This lead to a series of papers
 on the geometric properties of $M_j$ with uniformly bounded sectional curvature
 and $\vol(M_j)\to 0$ by Cheeger-Gromov, Fukaya, Rong, Shioya-Yamaguchi and others 
\cite{Fukaya-collapsing} \cite{CheegerGromovI} \cite{CheegerGromovII}
\cite{ShioyaYamaguchi-I} 
\cite{ShioyaYamaguchi-II}
 (c.f. \cite{Rong-collapsing-survey} \cite{Morgan-Tian-text}).  Collapsing
 Riemannian manifolds with boundary have been studied by Alexander-Bishop and
 Cao-Ge  \cite{Cao-Ge} 
 \cite{Alexander-Bishop-ThinI}
\cite{Alexander-Bishop-ThinII}.
This collapsing theory is an essential component of Perelman's
 proof of the Geometrization Conjecture using Hamilton's
 Ricci flow \cite{Morgan-Tian-text} \cite{Kleiner-Lott-perelman} \cite{Cao-Zhu-perelman}.

 In 1992,  Greene-Petersen applied the notion of Gromov's filling volume
 to find a lower bound on the volume of a ball in a Riemannian
 manifold with a uniform geometric contractibility function.   A geometric contractibility function is a function $\rho:(0,r_0] \to (0,\infty)$  with $\lim_{r\to 0}\rho(r)=0$ such that 
 any ball $B_p(r)\subset M$ is contractible in $B_p(\rho(r))\subset M$.
Applying Gromov's compactness theorem, one may conclude that
a sequence $M_j$ with a uniform geometric contractibility
function and a uniform upper bound on volume, has a subsequence that
converges in the Gromov-Hausdorff sense to a compact metric space.
Note that in this setting, volume is uniformly
bounded below so the sequence is not collapsing \cite{Greene-Petersen}.
The limits of such sequences of spaces have been studied by Ferry,
Ferry-Okun, Schul-Wenger and Sormani-Wenger \cite{Ferry}
\cite{Ferry-Okun}\cite{SorWen1} and \cite{SorWen2}.

Noncollapsing sequences of Riemannian manifolds with Ricci
curvature bounded from below were studied by Colding and
Cheeger-Colding in their work on almost rigidity in which they
weakened the conditions of well known rigidity theorems.  One example
of a rigidity theorem is the fact that 
 if $M^m$ has $\Ricci \ge (m-1)$ and $\vol(M^m)=\vol(S^m)$
 then $M^m$ is isometric to $S^m$.  Colding 
 proved a corresponding almost rigidity  theorem
 which states that if $M^m$ has $\Ricci \ge (m-1)$ then
 $\forall \epsilon>0$, $\exists \delta_{m,\epsilon}>0$ such that
 $\vol(M^m)\ge \vol(S^m)-\delta_{m,\epsilon}$ implies $d_{GH}(M^m, S^m)< \epsilon$.
 For a survey of such almost rigidity theorems see \cite{Colding-survey},
 \cite{ChCo-almost-rigidity} and Section 8 of \cite{Sor-length}.  
 The proofs generally involve an explicit construction of an almost
 isometry using distance functions and solutions to elliptic
 equations on the Riemannian manifolds.

When a sequence of $M_j$ with uniform lower bounds
on Ricci curvature does not collapse Cheeger-Colding proved
the volume is continuous and the Laplace spectrum converges.
Fukaya observed that when a sequence of Riemannian manifolds
with uniform lower bounds on Ricci curvature is collapsing, then
the Laplace spectrum does not converge.  In Figure~\ref{fig-Fukaya}
there are two sequences of Riemannian manifolds $A_i$ and
$B_i$ with nonnegative Ricci curvature converging to the same limit space
$A=B=[0,1]$ with the standard metric.  The eigenvalues converge
to real numbers, 
\be
\lambda^A_i:=\lim_{j\to\infty} \lambda_i(A_j) \textrm{ and }
\lambda^B_i:=\lim_{j\to\infty} \lambda_i(B_j)
\ee
but $\lambda^A_i \neq \lambda^B_i$.  While $\lambda^B_i$
made sense as eigenvalues for $[0,1]$, the other collection
of numbers $\lambda^A_i$ did not.  The alternating
sequence $\{A_1, B_1, A_2, B_2...\}$ also converges in the
Gromov-Hausdorff sense  but the eigenvalues do not converge
at all. \cite{Fukaya-87}.

This example disturbed Fukaya in light of the successful isospectral
compactness theorems of Osgood-Phillips-Sarnak,
Brooks-Perry-Petersen and Chang-Yang \cite{Osgood-Phillips-Sarnak} 
\cite{Brooks-Perry-Petersen} \cite{Chang-Yang-isospec-1990}
all of which imposed stronger conditions on the manifolds and
involved smooth convergence or even conformal convergence
of the sequences.  This lead to Fukaya's notion of metric measure
convergence.

\subsection{Metric Measure Convergence of Metric Measure Spaces} \label{Sect-2.2}
In 1987, Fukaya introduced the first notion of metric-measure convergence.
A sequence of metric measure spaces $(X_j, d_j, \mu_j)$
converge in the metric measure sense to a metric measure space 
$(X, d,\mu)$ if there is a sequence of $1/j$ almost isometries,
$f_j: X_j \to X$, such that push forwards of the measures, $f_{j*}\mu_j$,
converge weakly to $\mu$ on $X$.  Recall
that $\varphi_\#\mu(A):=\mu(\varphi^{-1}(A))$.

\begin{figure}[h] 
   \centering
   \includegraphics[width=4.6in]{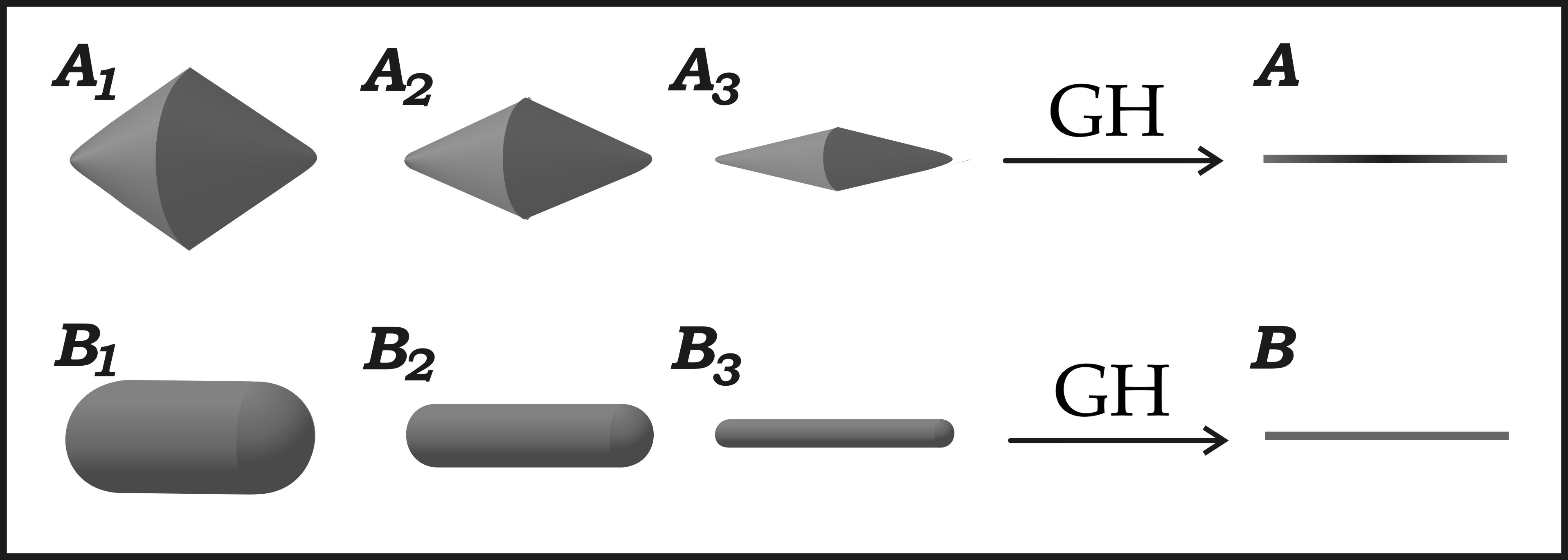} 
\caption{Metric Measure Convergence}
   \label{fig-Fukaya}
\end{figure}

In Figure~\ref{fig-Fukaya},  $A_j$ and $B_j$ are given
probability measures, $\mu_{A_j}$ and
$\mu_{B_j}$ proportional to $\mathcal{H}^2$.   Then
$B_j$ converge in the metric measure sense to $[0,1]$
with the standard metric and $\mu_B(W)=\mathcal{H}^1(W)$.  
Meanwhile $A_j$ converge in the metric measure sense to $[0,1]$ with the
measure 
\be
\mu_A(W)\,:=\, \int_W (2-|2-4x|)\,dx.
\ee
Defining the Laplacian
with respect to these measures, the spectrum for
$(A, \mu_A)$ is $\{\lambda^A_i\}$ and the spectrum for $(B, \mu_B)$
is $\{\lambda^B_i\}$.\footnote{See \cite{Fukaya-1987} or \cite{ChCo-PartIII} for
more details.} 

Cheeger-Colding then proved that any sequence of compact Riemannian
manifolds with uniform lower bounds on Ricci curvature endowed
with a probability measure proportional to the Hausdorff measure,
has a subsequence which converges in the metric measure sense.
The measures on the limit space satisfy the Bishop-Gromov comparison
theorem and are therefore doubling \cite{ChCo-PartI}.  
Sturm and Lott-Villani extended the notion
 of a Ricci curvature bound to general metric measure spaces using mass transport
 \cite{Sturm-06}\cite{Lott-Villani-09} (c.f. \cite{Villani-text}).
 Recently Topping and others have been developing a notion of Ricci flow on
 this larger class of spaces in hopes of defining Ricci flow through a singularity
 \cite{Topping-new}.
 
Cheeger-Colding also prove a Poincare inequality
on these limit spaces and deduced that the eigenvalues converge 
\cite{ChCo-Partii} \cite{ChCo-PartIII}.  They proved that the  
limit spaces of a sequence of Riemannian manifolds
with uniform lower bounds on Ricci curvature also have a notion
of tangent plane almost everywhere.  The tangent plane
at a point is found by rescaling the space $Y$ outward
and taking a pointed Gromov-Hausdorff limit.  At many points, one
does not get a unique limit under rescaling and the limits are
not necessarily planes or cones.  However they do exist at every point
and are called tangent cones.  Cheeger-Colding proved
tangent cones are Euclidean planes almost everywhere.  In fact
they prove $Y$ is a countably $\mathcal{H}^m$
rectifiable space and is a $C^{1,\alpha}$ manifold
away from the singularities. 
Key steps in the proof involve the Splitting
Theorem and the Poincare inequality (c.f. \cite{Cheeger-survey}).
There is an example where the limit space has infinite topological type
so it isn't a $C^{1,\alpha}$ manifold \cite{Menguy-inf-top-type}.
However, at least the manifold has a universal cover
unlike the Hawaii Ring depicted in Figure~\ref{fig-GHconv} \cite{SorWei1}.

It is also natural to study metric measure convergence without
Ricci curvature bounds on the sequence of manifolds.  As long as
the measure is doubling one can apply Gromov's compactness
theorem.  Sometimes the measure on the limit space is supported
on a smaller set.  This occurs for example for the sequence $B_j$
of Figure~\ref{fig-GHconv}.  Viewed as a metric measure limit
space, the limit space $B$ has a measure which is supported on the
two spheres.  

In 1981, Gromov introduced $\square_\lambda$ convergence to handle
this issue \cite{Gromov-metric}.  When $\lambda=1$, the limit of the sequence of $B_j$ is
just the two spheres with the line segment removed.  Basically, if $(Y,d,\mu)$
is the metric measure limit of a sequence, and if there is a ball such that
$\mu(B_y(r))=0$ then $y$ is removed from the space.  This leaves us
with a new, smaller limit space, $(\bar{Y}, d_Y, \mu)$, which may no longer
be connected.  This process may be called ``reduction of measure".
The spaces are then no longer close in the Gromov-Hausdorff topology.

In 2005, Sturm introduced a new distance between metric measure spaces
which leads to such limits more naturally and also interacts well with 
mass transport notions mentioned above \cite{Sturm-2006-I}.  He uses the
Wasserstein distance, $W_p$, to measure the distance between measures
in a set $Z$ and defines an intrinsic Wasserstein distance using an infimum
over all metric spaces $Z$ and all isometric embeddings $\varphi_i: X_i \to Z$:
\be \label{eqn-measure}
d_{W_p}((X_1,d_1, \mu_1), (X_2, \mu_2, d_2)) := 
\inf\{ d_{W_p}(\varphi_{1,\#}\mu_1, \varphi_{2\#}\mu_2)\}.
\ee
Sturm proved convergence with respect to this distance is equivalent to
Gromov's $\square_1$ convergence \cite{Sturm-2006-I}.  

Villani recently defined an intrinsic Prokhorov
distance replacing the Wasserstein distance of order $p$ in (\ref{eqn-measure})
with the Prokhorov
distance between measures.  As convergence with respect to the Wasserstein
distance and with respect to the Prokhorov distance  agree with weak convergence,
the intrinsic Prokhorov and intrinsic Wasserstein limits agree as well.  He refers
to these kinds of convergence as ``measure convergence"
Two metric measure spaces are a zero distance apart with respect to these
intrinsic measure distances  iff there is a measure preserving isometry between them
\cite{Villani-text}.

Villani also described  ``metric measure distances".  The Gromov-Hausdorff
Wasserstein distance is defined: $ d_{W_p}((X_1,d_1, \mu_1),(X_2, d_2, \mu_2)) :=$
\be
=\inf\{ d_H(\varphi_1(X_1), \varphi_2(X_2)) +
d_{W_p}(\varphi_{1\#}\mu_1, \varphi_{2\#}\mu_2)\},
\ee
where the infimum is taken over all metric spaces $Z$
and all isometric embeddings $\varphi_i:X_i\to Z$.
The Gromov-Hausdorff Prokhorov distance has the same formula replacing the
Wasserstein distance  by the Prokhorov distance between measures.  
Convergence with respect to these metric measure distances agrees with 
Fukaya's metric measure convergence. The advantage of
having multiple distances is that one may be easier estimate 
than another \cite{Villani-text}.

A sequence of metric measure spaces with doubling measures has a
converging subsequence and the
metric measure limit  agrees with the measure limit
of a sequence \cite{Villani-text}.  Any sequence of Riemannian manifolds
with lower bounds on their Ricci curvature satisfies the doubling condition
by the Bishop-Gromov volume comparison theorem.  The sequence of
$B_j$ in Figure~\ref{fig-GHconv} does not, as the metric measure and
measure limits do not agree.  

Despite the immense success in applying these definitions of convergence to
 study manifolds with Ricci curvature bounds, there has been a need to introduce
 a weaker form of convergence to study sequences of manifolds which do not
 satisfy these strong conditions.  Mathematicians studying manifolds with scalar
 curvature bounds and those interested only in sequences with an upper bound
 on volume and diameter without curvature bounds, need a weaker version of
 convergence.  In particular, geometric analysis related to cosmology and the
 study of the spacelike universe requires a weaker form of convergence.

An important example introduced by Ilmanen may be called 
the $3$-sphere of many splines.
It is a sequence of three dimensional spheres with positive scalar curvature whose
volume converges to the volume of the standard three sphere, but has an increasingly
dense set of thinner and thinner splines of ``length" $1$.  Cosmologically, one may
think of these splines as deep gravity wells.
The sequence does not converge in the
Gromov-Hausdorff sense because balls of radius $1/2$
centered on the tip of each spline are disjoint and the number of such disjoint
balls approaches infinity.  This example naturally lead Sormani and Wenger to develop
a notion of intrinsic flat convergence.

\subsection{Intrinsic Flat Convergence of Integral Current Spaces} \label{Sect-2.3} 
In 2008 Sormani and Wenger introduced the intrinsic flat distance
between compact oriented Riemannian manifolds \cite{SorWen1}\cite{SorWen2}:  
 \be
d_{\mathcal{F}}(M_1, M_2):=\inf \Big\{ d^Z_{\mathcal{F}}(\varphi_{1\#}T_1, \varphi_{2\#}T_2)
 \,{\bf :}
 \,\, isom\,\, \varphi_i:M_i \to Z\Big\}
 \ee
 where the flat distance in $Z$ is defined as in (\ref{eqn-def-flat}),
 the $\varphi_i:M_i\to Z_i$ are isometric embeddings as in (\ref{eqn-isom-embed})
 and where $T_i$ are defined by integration over $M_i$ so that
 $\varphi_{i\#}T_i (f, x_1,...x_k) = $
 \be \label{eqn-push-forward}
 =
 T_i (f\circ \varphi_i, x_1\circ\varphi_i ,... x_k\circ\varphi_i)
 = \int_{M_i} f\circ \varphi_i \,\, d(x_1\circ \varphi_i) \wedge \cdots \wedge d(x_k\circ \varphi_i).
 \ee
 
 One may immediately note that 
 $d_{\mathcal{F}}(M_1, M_2) \le \vol(M_1)+\vol(M_2)<\infty$
 as we may always take the integral current $A$ in (\ref{eqn-def-flat})
 to be $\varphi_{1\#}T_1-\varphi_{2\#}T_2$ and $B=0$.  
 The intrinsic flat distance is a distance between compact oriented
 Riemannian manifolds in the sense that $d_{\mathcal{F}}(M_1,M_2)=0$
 iff there is an orientation preserving isometry between $M_1$ and $M_2$ \cite{SorWen2}.  
  
 Note that in practice it is often possible to estimate the intrinsic flat distance
 using only notions from Riemannian geometry.  That is, if $M^k_1$ and $M^k_2$
 isometrically embed into a Riemannian manifold $N^{k+1}$ such that 
 $
 \partial N^{k+1}=\varphi_1(M^k_1) \cup \varphi_2(M^k_2)
 $
 and the manifolds have been given an orientation consistent
 with Stoke's theorem on $N^{k+1}$, then
 $
 d_{\mathcal{F}}(M_1, M_2) \le \vol(N^{k+1})
 $
 This viewpoint makes it quite easy to see that Ilmanen's $3$ sphere of many splines
 example describes a sequence converging to the standard sphere.  Note that
 one cannot just embed the sequence into four dimensional Euclidean space and
 take $N^{4}$ to be the flat region lying between the two spheres because such
 $\varphi_i$ would not be isometric embeddings.  
 Instead, one rotates each spline into half a thin $4$ dimensional
 spline and glues it smoothly to a short $S^3 \times [0,\epsilon_j]$ where 
 $\epsilon_j$ is ``thinness" of the spline.  This produces a four dimensional
 manifold $N^4$ with $S^3$ isometrically embedded as one boundary and
 the $3$ sphere of many splines isometrically embedded as the other boundary.
 So the intrinsic flat distance between these two spaces is $\le \vol(N^4)$
 which is approximately $\epsilon_j \vol(S^3)$ plus the sum of the volumes
 of the four dimensional splines each of which is approximately
 $ \epsilon_j^3$.  See \cite{SorWen2} for this example and many others.

 
 The limit spaces are called integral
 current spaces, and are oriented weighted countably $\mathcal{H}^k$
 rectifiable metric spaces.  They have the same dimension as the
 sequence.  They have a set of
Lipschitz charts as well as a notion of an approximate tangent space.
The approximate tangent space is a normed space whose
norm is defined by the metric differential.  The notion of
the metric differential was developed in work of 
Korevaar-Schoen \cite{Korevaar-Schoen} and Kirchheim \cite{Kirchheim}.  
 See \cite{SorWen2} for
 a number of examples of limit spaces.


 Integral current spaces are said to have a ``current structure", $T$,
 defined by integration using the orientation and weight.  One writes
 $(X,d,T)$.  The current structure defines a mass measure $||T||$.
 The points in $x$ all have positive density with respect to this measure.
 These spaces have a notion of boundary coming
 from this current structure.  The boundary is also an integral current space.    
  \footnote{See \cite{SorWen2} for more details about the relationship
  between $X$ and $T$ and how to find the boundary.}
  
    An important integral current space is the $0$ current
 space, and collapsing sequences of manifolds,
 $\vol(M_j)\to 0$, converge to the $\bf{0}$ space.
 Notice that in Figure~\ref{fig-GHconv}, the sequences $A_j$, $C_j$
 and $D_j$ all converge to $0$ in the intrinsic flat sense.  
 The sequence $B_j$ converges to a pair of spheres in the intrinsic
 flat sense where the limit does not include the line segment.


Most results about integral current spaces and intrinsic flat convergence
are proven by finding a way to isometrically embed the entire sequence
into a common metric space, $Z$, and apply the theorems for integral
currents in $Z$.  The rectifiablility, dimension and boundary properties of
the limits follows from Ambrosio-Kirchheim's corresponding results on
integral currents.  So does the slicing theorem and the lower semicontinuity 
of mass.  Combining Gromov's compactness theorem with Ambrosio-Kirchheim's
compactness theorem, one sees that a sequence of manifolds which converge
in the Gromov-Hausdorff sense that have a uniform upper bound on volume,
also converge in the intrinsic flat sense \cite{SorWen2}.

In general the intrinsic flat and Gromov-Hausdorff limits
do not agree: the intrinsic flat limit may be a strict subset of the Gromov-Hausdorff 
limit.  If the Gromov Hausdorff limit is lower dimensional than the sequence,
then the intrinsic flat limit is the ${\bf{0}}$ space.  Any regions in the Gromov-Hausdorff
limit that are lower dimensional disappear from the intrinsic
flat limit  \cite{SorWen2}.
Note that this is in contrast with
metric measure limits, which do not lose regions of lower dimension that
still have positive limit measures.  Such a situation can occur, for example, if
a thin region collapsing to lower dimension is very bumpy and has a uniform
lower bound on volume.

Sequences of Riemannian manifolds may also disappear due to cancellation.
Recall that in Euclidean space, a submanifold which curves in on itself as in
Figure~\ref{fig-cancels} will have cancellation in the limit due to the opposing
orientation.  An example described in \cite{SorWen1} is a sequence of
three dimensional manifolds with positive scalar curvature created by taking a
pair of standard three spheres and joining them by increasingly dense and
increasingly thin and short tunnels.  The Gromov-Hausdorff limit is a 3-sphere
while the intrinsic flat limit is the $0$ current space.  If each tunnel is twisted,
then the intrinsic flat limit is the standard three sphere with multiplicity two.

In some
cases the intrinsic flat limits and Gromov-Hausdorff limits agree giving
new insight into the rectifiability of the Gromov-Hausdorff limits. 
They agree for 
noncollapsing sequences of manifolds with nonnegative Ricci curvature and also for
sequences of manifolds with uniform linear contractibility 
functions and uniform upper bounds on volume \cite{SorWen1}.  As a consequence the
Gromov-Hausdorff limits are countably $\mathcal{H}^m$ rectifiable metric spaces.
This was already shown by Cheeger-Colding for the sequences with
bounded Ricci curvature  but is
a new result for the sequences with the linear geometric contractibility hypothesis.
With only uniform geometric contractibility functions that are not linear,
Schul and Wenger have shown the limits need not be so rectifiable \cite{SorWen1}
and when there is no upper bound on volume, Ferry has shown the limit 
spaces need not even be finite dimensional \cite{Ferry}.  
Recall that Cheeger-Colding have proven limits of manifolds with nonnegative
Ricci curvature have Euclidean tangent cones almost everywhere.  In contrast,
there is a sequence of Riemannian manifolds with uniform
linear contractibility functions that converge to the taxicab space \cite{SorWen2}.

The key ingrediant in the proofs of these noncancellation results is an estimate
on the filling volumes of small spheres in the spaces.  The filling volumes of
the spheres are continuous with respect to the intrinsic flat distance and can
be more easily applied to control sequence than the mass (or volume).
It is conjectured that a sequence of three dimensional Riemannian manifolds
with positive scalar curvature, no interior minimal surfaces and a uniform upper
bound on volume and on the area of the boundary converges in the intrinsic
flat sense without cancellation.  Such manifolds are important in the study of
general relativity. \cite{SorWen1}

In fact, Wenger has proven that a sequence of
oriented Riemannian manifolds with a uniform upper bound on diameter and
on volume has a subsequence which converges in the intrinsic flat sense
to an integral current space \cite{Wenger-compactness}.  The proof involves
an even weaker notion of convergence: the Ultralimit.

\subsection{Ultralimits} \label{Sect-2.4}

The 
notion of an ultralimit was introduced by van den Dries and Wilkie in 1984
\cite{DW-ultralimit}
and developed by Gromov in \cite{Gromov-asympt}.  An ultralimit of a sequence
of metric spaces is defined using Cartan's 1937 notion of a nonprincipal ultrafilter: 
a finitely
additive probability measure $\omega$ such that all subsets $S \subset \NN$
are $\omega$ measureable, $\omega(S)\in \{0,1\}$ and
$\omega(S)=0$ whenever $S$ is finite.  Given an ultrafilter, $\omega$, and a bounded
sequence of $a_j \in \RR$, there exists an ultralimit, $L=\omega-\lim_{j\to\infty} a_j$,
such that $\forall \epsilon >0$, $\omega\{j: |a_j-L|<\epsilon \}=1$.  

Given a sequence of compact metric spaces, $(X_j, d_j)$, with a uniform bound
on diameter, and given an ultrafilter, $\omega$, there is an ultralimit $(X,d)$
which is a metric space that may no longer be compact, but is at least complete.
The space, $X$, 
is constructed as equivalence classes of sequences $\{x_j\}$ where $x_j \in X_j$
and the metric of $X$ is defined by taking ultralimits:
\be
d(\{x_j\}, \{y_j\}) := \omega - \lim_{j\to\infty} d_j(x_j,y_j).
\ee
Two sequences are equivalent when the distance between them is zero.  Notice
that one does not need a subsequence to find a limit and that in general the
ultralimit depends on $\omega$.  

If $X_j$ has a Gromov-Hausdorff limit, then the ultralimit is the Gromov-Hausdorff
limit and there is no dependence on $\omega$.  If one has two 
sequences $X_j \GHto X$ and $Y_j\GHto Y$, then the ultralimit
of the alternating sequence $\{X_1, Y_1, X_2, Y_2, X_3, ...\}$ is $Y$ iff
$\exists N $ such that $\forall n>N$ we have
$\omega\{2n, 2n+2, 2n+4, 2n+6,...\}=1$.  Otherwise the ultralimit is $X$.  

It is often useful to compute the ultralimits of sequences
which do not have Gromov-Hausdorff limits.   For example,
the three sphere of many splines sequence, $\{M_j\}$,
converges to a standard three sphere
with countably many unit line segments attached at various points on the sphere.
To see this, one may view each $M_j$ as a union of regions $W_j \cup U_{1,j}\cup U_{2,j}\cup...\cup U_{N_j,j}$, 
where each $U_{i,j}$ covers a spline and $W_j$ covers the spherical
portion.  Then $W_j \GHto S^3$, while $U_{i_j,j} \GHto [0,1]$.  Thus any
ultralimit of these $M_j$ is built by connecting countably many intervals to a three
sphere.     The locations where
the intervals are attached to the sphere may depend on the ultrafilter, $\omega$.

If $X_j$ are geodesic spaces, 
then the ultralimit is a geodesic space.
If the $X_j$ are $CAT(0)$ spaces then the ultralimit is as well.
Ultralimits have been applied extensively in the study of CAT(0) spaces
and also Lie Groups.  See for example \cite{BH} and \cite{Kapovich-text} for
more details as well as applications.  

\section{Speculation}
While the methods of intrinsic convergence defined above have all 
proven to be useful in a variety of settings, and should in fact have
more applications that have not yet been discovered or fully explored,
it is clear that each has its disadvantages.  Gromov-Hausdorff
convergence, like Hausdorff convergence, provides no rectifiability
for its limits.  The intrinsic flat convergence, like Federer-Fleming's flat
convergence, has difficulty with cancellation in the limit.  Ultraconvergence
has the same difficulties as Gromov-Hausdorff convergence and the
additional concern that the spaces in the sequence are not particularly
close to the limit space in a measurable way using a notion of distance.

There are important problems both related to General Relativity
and Ricci flow which have not yet been addressed using the above methods.
While intrinsic flat convergence may prove useful as it is further explored,
other methods of convergence may be necessary to address all the problems
that arise.


\subsection{Intrinsic Varifold Convergence}\label{Sect-3.1}
For some time, mathematicians have been exploring possible methods of extending Ricci
flow through singularities.  Mean curvature flow behaves a lot like Ricci flow
and the Brakke flow uses varifolds to extend this notion to the nonsmooth
setting. The key missing ingredient that has kept people from directly extending the work
on Brakke flow to better understand Ricci flow is that there has never been
a notion of an intrinsic varifold space and intrinsic varifold convergence.

In the prior section we saw how mathematicians have repeatedly applied
Gromov's trick of isometrically embedding metric spaces into a common space,
measuring the distance between them in that common space, and then taking
an infimum.  Here, however,
there is no metric measuring the distance between varifolds, just the
natural notion of weak convergence against test functions.  
One might try to define a notion of distances between varifolds
in Euclidean space which defines a convergence that agrees with the standard 
varifold convergence and then apply Gromov's trick.  

Alternately, one might venture to say that a sequence of Riemannian manifolds, $M_j$,
converges to a space, $M$, as varifolds if and only if there is a sequence of isometric
embeddings from, $M_j$, into a common space, $Z$, such that the images of
$M_j$ converge as varifolds in $Z$ to an image of $M$.  However, there is as yet
no notion of convergence as varifolds in a metric space.  Recall that if
$Z$ is Euclidean space, convergence as varifolds, means weak convergence
as measures on $\RR^N \times \Gamma(N,k)$.  We would need some sort
of notion of a $\Gamma(Z,k)$ perhaps representing all possible ``tangent spaces"
of $Z$ at a point.  Perhaps this might be done by first isometrically embedding
$Z$ into a Banach space.

\subsection{Area Convergence} \label{Sect-3.3}

One may view Gromov-Hausdorff convergence as a notion of convergence
defined by the fact that lengths of minimizing geodesics converge.  
Measure convergence is defined by the fact that volumes converge.
It becomes natural to wonder, particularly in the case of three dimensional
manifolds, whether there is a notion of convergence defined by the fact that
areas, or areas of minimal surfaces, converge.

Note that Ilmanen's example of the three sphere of many splines, $\{M_j\}$,
converges in some area sense to the standard three sphere.  There are
radial projections $f_j: M_j \to S^3$ which are almost area preserving, 
one-to-one and onto.  Such a map seems perhaps to extend the notion
of an almost isometry without involving a distance.  However requiring
a diffeomorphism to define area convergence seems too strong for any
applications.   

Burago, Ivanov and Sormani spent a few years investigating possible notions
of area convergence and the area distance between spaces.   They
observed that two distinct metric spaces could easily have a surjection between
them that is area preserving.  One space could, for example, be the unit
square, $X=[0,1]\times[0,1]$, with the Euclidean metric.  The other space
could be the unit square with a ``pulled thread": $Y=[0,1]\times[0,1]/ \sim$ where 
$(x_1,x_2)\sim (y_1,y_2)$ iff $x_1=y_1=0$ with the metric:
\be
d_Y((x_1,x_2), (y_1,y_2)) =\min \{ \sqrt{(x_1-y_1)^2+(x_2-y_2)^2}, |x_1|+|y_1| \}.
\ee
This seemed to indicate the need for a notion of ``area space" defined as
an equivalence class of metric spaces.  
It was essential before even beginning the project
to verify that each equivalence class would only include one Riemannian 
manifold.

Another difficulty arose in that one would not expect to always find an
almost area preserving surjection between spaces which should intuitively
be close in the area sense.  In an almost isometry, there are two requirements:
almost distance preserving and almost onto.  Both requirements involve distance
and one can prove that if there is an almost isometry from a metric space
$X$ to a metric space $Y$ then one can find an almost isometry in the reverse
direction (although the errors change slightly in reverse) \cite{Gromov-metric}.
To take advantage of Gromov's ideas, it was decided that one needs to convert
the notion of area into a notion of distance.  

Given a Riemannian manifold, $M$, one may examine the loop
space, $\Omega(M)$.  The flat distance between loops is defined using
the notion of area.  Recall that the flat distance defined in (\ref{eqn-def-flat})
involves both area and length, however, in the setting of closed
loops one may choose $A=0$ and only take an infimum over all
possible fillings by surfaces $B$.  Perhaps an area convergence of 
$M_j$ to some sort of area space $M$ could be defined by taking
pointed Gromov-Hausdorff limits or  ultralimits of the $\Omega(M_j)$.

Preliminary work in this direction was completed by Burago-Ivanov.  They
proved that if a pair of three dimensional Riemannian manifolds, $M^3$
and $N^3$, have isometric loop spaces, then $M^3$ and $N^3$ are isometric
\cite{Burago-Ivanov-area}.  This step alone was very difficult because $M^3$
is not isometrically embedded into its loop space, $\Omega(M)$, so one must
localize points and planes in $TM$ using sequences of loops without controlling
the lengths of the loops.    

Now one may move forward and investigate
whether there is a natural compactness theorem for some notion
of area convergence.   It would not be surprising if
positive scalar curvature on $M_j$ control areas of minimal
surfaces well enough to control $\Omega(M_j)$.  Or perhaps one
could skip compactness theorems, and apply ultralimits to the
$\Omega(M_j)$ and see what properties are conserved under
such a limit.  Ultimately one does not just want to take a limit of
$\Omega(M_j)$ but find an area space $X$, such that $\Omega(X)$
is the limit of the $\Omega(M_j)$.  Then one could say $X$ is
the area limit of $M_j$.

\subsection{Convergence of Lorentzian manifolds} \label{Sect-3.2}
Another direction of research that is fundamental to General Relativity
is the development of weak forms of convergence for Lorentzian manifolds.
While Gromov-Hausdorff convergence has proven useful in the
study of the stability of the spacelike universe \cite{Sor-Cosmos}, one needs
to extend Gromov-Hausdorff convergence to the Lorentzian setting
to study the stability of the space-time universe.
The techniques described above immediately fail for these spaces
because they do not isometrically embed into metric spaces.  

Noldus has developed a Gromov-Hausdorff distance between Lorentzian
spaces \cite{Noldus-Lorentzian}.
Further work in this direction including a notion of the arising limit
spaces appears in \cite{Noldus-Limit}
and compactness theorems and an examination of causality in the limit spaces are
proven with Bombelli in \cite{Bombelli-Noldus}.   This work has not yet been applied to
prove stability questions
arising in general relativity.

Perhaps similar methods might be applied to extend the notion of the intrinsic
flat distance to Lorentzian spaces.  

\subsection{Acknowledgements}
The author is indebted to Jeff Cheeger for all the excellent advise he has given
over the past fifteen years.  His intuitive explanations of the beautiful geometry
of converging Riemannian manifolds in his many excellent talks have been
invaluable to the mathematics community.   The author would also like to thank
the editors for soliciting this article and William Wylie for discussions regarding the exposition.

\bibliographystyle{plain}
\bibliography{2010}
\end{document}